\documentclass[onecolumn,showpacs,preprintnumbers,amsmath,amssymb]
{revtex4-1}
\usepackage{epsfig}
\usepackage{graphicx,epstopdf}
\usepackage{placeins}

\usepackage{amssymb}
\usepackage{amsthm}

\usepackage{amsmath,amssymb,amsfonts,latexsym,cancel}
\usepackage{rawfonts}
\usepackage{pictexwd}
\usepackage[all]{xy}
\newtheorem {theorem}{Theorem}
\newtheorem {lemma}{Lemma}

\begin{document}


\title{The Moser's formula for the division of the circle by chords problem revisited}

\author{Carlos Rodr\'{\i}guez-Lucatero}
\email{crodriguez@correo.cua.uam.mx}
\affiliation{Departamento de Tecnolog\'{\i}as de la Informaci\'on, Universidad Aut\'onoma Metropolitana-Cuajimalpa,
Torre III,
Av. Vasco de Quiroga 4871,
Col.Santa Fe Cuajimalpa, M\'exico, D. F.,
C.P. 05348, M\'exico}

\begin{abstract}
The enumeration of the regions formed when circle  is divided by secants drawn from points on the circle is one of the examples where the inductive reasoning fails as was pointed out by Leo Moser in the Mathematical Miscellany in 1949. The formula that gives the right number of regions can be deduced by combinatorics reasoning using the Euler's planar graph formula, etc. My contribution in the present work is to reformulate and solve such problem in terms of a  fourth order difference equation and to obtain the formula proposed by Leo Moser. 
\\
\\
\\
{\em Mathematics Subjects Classification}: 05A15
\\
\end{abstract}

\maketitle


\section{Introduction}
A problem sometimes known as Moser's circle problem asks to determine the number of pieces into which a circle is divided if $m$ points on its circumference are joined by chords with no three internally concurrent. 
The number of regions formed inside the circle  when it is divided by the chords as mentioned, can be sketched for the fisrt $5$ steps in the following figure:

\begin{figure}[htbp]
\begin{center}
\includegraphics[width=400pt]{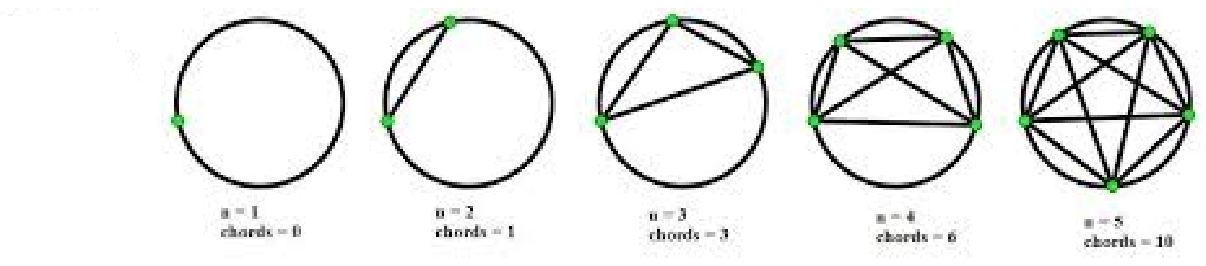}
\end{center}
\end{figure}

If we label the regions formed by this division by secants, the sequence of the number of regions generated  in terms of number of points till this point is $1,2,4,8,16$. If the number of points is denoted by $m$ and we try to guess the functional behavior of the sequence the induction tell us that it is $2^{m-1}$. By following this procedure if we add one more point and draw the corresponding secants we get $31$ regions instead of $32$. Let me tabulate this behavior for the seven first points in the following table

\begin{tabular}{|c|c|}
 \hline 
 Points  & regions \\
 $m$     & $f(m)$ \\
 \hline 
 1 & 1 \\ 
 \hline 
 2 & 2 \\ 
 \hline 
 3 & 4 \\ 
 \hline 
 4 & 8 \\ 
 \hline 
 5 & 16 \\ 
 \hline 
 6 & 31 \\ 
 \hline 
 7 & 57 \\ 
 \hline  
 \end{tabular}

As can be seen the functon $2^{m-1}$ no longer describes the behavior of the sequence. That is the reason why Leo Moser in \cite{Moser1} pointed that the inductive method for guessing the next element in a numerical sequence can fail.  In fact the title of the section in \cite{Moser1} was {\em On the danger of induction}. He leaves to the reader as an exercise, to show that the number of regions formed by joining points on a circle by chords is $f(m)=\sum_{j=0}^{m}\binom{m-1}{j}$.

The function that describes the behavior of the sequence is:

\begin{equation}
f(m)=\frac{m^4-6m^3+23 m^2-18m+24}{24}
\end{equation}

This can be proven in many different forms. In the following sections I will describe how it was demonstrated
 in \cite{Maier1} using combinatorial arguments, in \cite{Jobbings1} by using the Euler's planar graphs formula, and finally how it can be proven by using a fourth degree difference equation.

\section{Deduction using a combinatorial argumentation}\label{combingeom}

In this section we will describe a proof of the Moser's formula based on the article \cite{Maier1}(see also\cite{Jobbings1}). 
In order to find the actual formula on the number of regions formed, the author of \cite{Jobbings1} states the following 
\\
{\bf Result:} The number of regions formed is

\begin{equation}\label{combinat1}
1+\binom{m}{2}+\binom{m}{4}
\end{equation}

From the combinatorics it is known that $\binom{r}{s}=\frac{r!}{s!(r-s)!}$ and if we apply that to \ref{combinat1} we obtain the Moser's formula

\begin{equation} \label{combinat2}
f(m)=\frac{m^4-6m^3+23m^2-18m+24}{24}
\end{equation}

The proof of \ref{combinat2} given in \cite{Jobbings1} uses the fact that the binomial coeficient $\binom{r}{s}$ counts the ways we can choose $s$ elements from a set of $r$ different elements.

In order to formalize the demonstration, I will start by some basic lemas that can be used for the proof.

\begin{lemma} \label{lemmacombinat1}
The total number of chords that can be created from the $m$ points on the circle is
\begin{equation}
\binom{m}{2}
\end{equation}
\begin{proof}
If we have $m$ points each chord is formed by taking two points of the $m$ lying on the perimeter of the circle. In other words the number of chords is determined by the number of ways in which the $m$ points on the circle can be taken to form a chord.
Then the total number of posible chords equals the total number of different ways in which  $m$ points can be taken in groups of $2$ what is equal to $\binom{m}{2}$.
\end{proof}
\end{lemma}

\begin{lemma} \label{lemmacombinat2}

The total number of interior crossing points of the chords obtained from the $m$ points on the circle is
\begin{equation}
\binom{m}{4}
\end{equation}
\begin{proof}
For calculating the number of interior intersection points inside the circle, we must count the number of ways in which the $m$ points in the circle can be taken such that their related chords intersect, which can be done taking four points to form their chords and their corresponding intersections, what gives a total of $\binom{m}{4}$ ways
\end{proof}
\end{lemma}

\begin{theorem} \label{theoremcombinat}
Let $m$ be the number of points on the circle and $f(m)$ the number of regions formed by the division of a circle by chords. The total number of regions formed is

\begin{equation} 
f(m)=\frac{m^4-6m^3+23m^2-18m+24}{24}
\end{equation}

\begin{proof}
Let me start by counting the regions formed by taking each chord one by one. Each time a new chord is drawn, it crosses a number of regions dividing them into two. The number of new created regions is equal to the number of regions crossed by the new chord, which is one more than the number of chords crossed. Given that the new chord cannot pass through a previously drawn point of intersection, then the number of chords crossed is equal to the number of interior points of intersection of the new chord. As consequence, the number of new regions created by this new chord equals the number of interior points of intersection of this chord plus one.
Because of that, taking into accont lemma \ref{lemmacombinat1} and lemma \ref{lemmacombinat2} the total number of new regions created by drawing all the chords is equal to the number of added chords plus the number of interior points of intersection as well as the fact that at the begining we have one region, the total number of regions formed is

\begin{equation} \label{combinat3}
f(m)=1+\binom{m}{2}+\binom{m}{4}=\frac{m^4-6m^3+23m^2-18m+24}{24}
\end{equation}

\end{proof}
\end{theorem}

\section{Deduction using the planar graph Euler's formula}

I found in \cite{Jobbings1} an elegant proof based on the famous planar graph Euler's formula that I will develop in this section. Let me start with the statement of the Euler's result
\newline
\\
\begin{theorem} \label{eulerformula}
 Let $V$ the number of vertices, $E$ the number of edges and $F$ the number of faces of a planar graph. Then $V-E+F=2$.
\end{theorem}

In order to use the planar graphs Euler's for deducing the Leo Moser's formula the author of \cite{Jobbings1} relate the set $V$ with the orginal points on the circle as well as with the interior intersection points, the set $E$ with the chords and arc's formed by the points on the circle and the set $F$ with the formed regions. 
We can formalize this method as follows

\begin{theorem} \label{theoremeuler}

Let $P$ be the set of $m$ points on the circle, $I$ the set of $\binom{m}{4}$ interior crossing points and $f(m)$ the number of regions formed by the division of the circle by chords. 
Let $G=(V,E)$ be the planar graph obtained from the circle division by chords, where $V=P \cup I$, $E$ the edges of the planar graph and $F=f(m)$. The number of regions or faces formed is

\begin{equation} 
f(m)=\frac{m^4-6m^3+23m^2-18m+24}{24}
\end{equation}
 
\begin{proof}
There are $m$ points on the circle and there are $\binom{m}{4}$ intersections of the chords. Then we have a total of
\begin{equation}\label{euler1}
V=(P \cup I)=m+\binom{m}{4}
\end{equation}

vertices. In order to count the total number of edges it must be noticed that we have $m$ circular arcs. Given that we have $\binom{m}{4}$ interior intresection points where four edges meet then we have $4\binom{m}{4}$ additional edges. Due to the fact that we have $\binom{m}{2}$ chords corresponding to two edges that meet the circle we have $2\binom{m}{2}$ more edges. Then we have $2\binom{m}{2}+4\binom{m}{4}$ edges generated by the chords, but due to counting process we have counted them twice. Then this quantity must be divided by $2$ giving a total of $\binom{m}{2}+2\binom{m}{4}$. Hence the total number of edges is

\begin{equation} \label{euler2}
E=m+\binom{m}{2}+2\binom{m}{4}
\end{equation}
The number of faces is related with the number of regions $f(m)$ as follows

\begin{equation} \label{euler3}
F=f(m)+1
\end{equation}

From the \ref{eulerformula} we know that

\begin{equation}\label{euler4}
V-E+F=2
\end{equation}

Replacing \ref{euler1},\ref{euler2},\ref{euler3} in \ref{euler4} we get

\begin{equation}\label{euler5}
\Bigg\{ m+\binom{m}{4} \Bigg\}-\Bigg\{ m+\binom{m}{2}+2\binom{m}{4} \Bigg\}+ f(m)+1 = 2
\end{equation}

Simplifying \ref{euler5} we get

\begin{equation}\label{euler6}
-\binom{m}{2}-\binom{m}{4} + f(m) = 1
\end{equation}

From \ref{euler6} we express $f(m)$ in terms of the other elements of the expression as follows

\begin{equation}\label{euler7}
f(m) = 1 +\binom{m}{2}+\binom{m}{4}=\frac{m^4-6m^3+23m^2-18m+24}{24} 
\end{equation}
 
and in that way we have finally mathematically proven the Moser's formula. 
\end{proof}
\end{theorem}

\section{My deduction of the Leo Moser's formula by solving a difference equation} \label{secgenfunc}

An alternative method that I propose for solving the {\em Moser's circle division by chords} is by obtaining a recurrence from the numerical sequence of the number of regions formed and then solve the related difference equation. For this end we have to obtain a recurrence relation to be solved. This recurrence relation can be obtained from the succesion of regions, using the technique of successive differences frequently applied in problems of inductive reasoning on numerical sequences. Applying the successive differences technique on the sequence $\{a_{0},a_{1},a_{2},a_{3},a_{4},a_{5}\}$ of six elements to obtain the seventh element $a_{6}$. Let me state this first result as a lemma.

\begin{lemma} \label{lemmarecurrence}
The elements of the sequence $1,2,4,8,16,31,57, \ldots $ that represent the growth behavior on the number of regions formed by dividing the circle by chords can be generated by a fourth degree recurrence relation
\begin{proof}
We start by applying the successive differences method for guessing the next element in the sequence on the first six elements of such numerical sequence.
\\
\begin{tabular}{|c|c|c|c|c|c|c|c|c|c|c|c|c|}
\hline 
$a_{0}$ &  &$a_{1}$& • &$a_{2}$& • &$a_{3}$& • & $a_{4}$ & • & $a_{5}$ & • &$a_{6}$\\
 
\hline 
1 &  & 2 & • & 4 & • & 8 & • & 16 & • & 31 & • & {\bf 57}\\
 
\hline 
• & 1 & • & 2 & • & 4 & • & 8 & • & 15 & • & • & • \\ 
\hline 
• & • & 1 & • & 2 & • & 4 & • & 7 & • & • & • & •  \\ 
\hline 
• & • & • & 1 & • & 2 & • & 3 & • & • & • & • & • \\ 
\hline 
• & • & • & • & 1 & • & 1 & • & • & • & • & • & • \\ 
\hline 
\end{tabular} 
\\
\\
From the calculations it can be noticed that we stop the successive differences procedure when the differences become constant. Then if we sum the last element of each row plus the last element on the original sequence of numbers we can obtain the next element on the sequence. So in the example, from the summation $1+3+7+15+31$ we obtain $57$ that correspond to the next value on the sequence.   
From the successive differences table we obtain the following recurrence relation

\small{
\begin{equation} \label{rec1}
\begin{array}{c}
a_{n+4} = a_{n+3} +(a_{n+3}-a_{n+2})+((a_{n+3}-a_{n+2})-(a_{n+2}-a_{n+1}))+ \\
(((a_{n+3}-a_{n+2})-(a_{n+2}-a_{n+1}))-((a_{n+2}-a_{n+1})-(a_{n+1}-a_{n})))+1
\end{array} 
\end{equation}
}

The expression \ref{rec1} is the desired fourth degree recurrence relationship.

\end{proof}
\end{lemma}

Simplifying and reordering the terms of the equation \ref{rec1} we get the next expression
\begin{equation} \label{rec2}
a_{n+4}-4 a_{n+3} + 6 a_{n+2} - 4 a_{n+1} + a_{n} = 1
\end{equation}

adding the inicial conditions to \ref{rec2} we get the following difference equation

\begin{equation} \label{rec3}
a_{n+4}-4 a_{n+3} + 6 a_{n+2} - 4 a_{n+1} + a_{n} = 1, ~~~~ a_{0}=1,a_{1}=2,a_{2}=4,a_{3}=8
\end{equation}

Once we have obtained the difference equation \ref{rec3} we can solve it by many existing methods \cite{Grimaldi1},\cite{Sedgewick1}. The purpose of solving \ref{rec3} is to releat it with the deduction of the Leo Moser's formula of the number of regions formed by the division of the circle by chords. In what follows I will do it in two different ways
\begin{itemize}
  \item By solving \ref{rec3} by the generating functions method
  \item By solving \ref{rec3} by the solution of non-homogeneous linear with constant coefficients difference equation method.
\end{itemize}

To formalize these results I will state them as the following theorems.

\begin{theorem}\label{theoremgenfunc}
Let $n$ the subindex of the recurrence \ref{rec3}, $m$ the number of points on the circle and $f(m)$ the number of regions formed by the division of the circle by chords. Let $n=m-1$. If we solve \ref{rec3} by the generating functions method the number of regions formed by the division of the circle by chords
is
\begin{equation}
f(m) = \frac{m^4-6m^3+23m^2-18m+24}{24} 
\end{equation}

\begin{proof}
From \ref{lemmarecurrence} we have \ref{rec3}.
The application of equation \ref{rec3} for different values of $n$ gives the following results
\\
\begin{equation} \label{rec4}
\begin{array}{ll}
(n=0)  & a_{4}-4 a_{3} + 6 a_{2} - 4 a_{1} + a_{0} = 1 \\
(n=1)  & a_{5}-4 a_{4} + 6 a_{3} - 4 a_{2} + a_{1} = 1 \\
(n=2)  & a_{6}-4 a_{5} + 6 a_{4} - 4 a_{3} + a_{2} = 1 \\
(n=3)  & a_{7}-4 a_{6} + 6 a_{5} - 4 a_{4} + a_{3} = 1 \\
\vdots & \vdots
\end{array}
\end{equation}

If we multiply \ref{rec4} by $x^{0}$ the first row, the second row by $x^{1}$, the third row by $x^{2}$, the fourth row by $x^{3}$ and so on we get
\\
\begin{equation} \label{rec5}
\begin{array}{ll}
(n=0)  & a_{4} x^{0}-4 a_{3} x^{0} + 6 a_{2} x^{0} - 4 a_{1} x^{0} + a_{0} x^{0} = x^{0} \\
(n=1)  & a_{5} x^{1}-4 a_{4} x^{1} + 6 a_{3} x^{1} - 4 a_{2} x^{1} + a_{1} x^{1} = x^{1} \\
(n=2)  & a_{6} x^{2}-4 a_{5} x^{2} + 6 a_{4} x^{2} - 4 a_{3} x^{2} + a_{2} x^{2} =  x^{2} \\
(n=3)  & a_{7} x^{3}-4 a_{6} x^{3} + 6 a_{5} x^{3} - 4 a_{4} x^{3} + a_{3} x^{3} = x^{3} \\
\vdots & \vdots
\end{array}
\end{equation}
\\
Summing up the rows of \ref{rec5} we obtain
\\
\begin{equation} \label{rec6}
 \sum_{n=0}^{\infty}a_{n+4} x^{n} - 4 \sum_{n=0}^{\infty}a_{n+3} x^{n} + 6 \sum_{n=0}^{\infty}a_{n+2} x^{n} - 4 \sum_{n=0}^{\infty}a_{n+1} x^{n} + \sum_{n=0}^{\infty}a_{n} x^{n} = \sum_{n=0}^{\infty} x^{n}
\end{equation} 
\\ 
Trying to equate the subindex of the coeficientes and the powers of the variables we rewrite \ref{rec6} as
\\
\begin{equation} \label{rec7}
\begin{array}{l}
 x^{-4} \sum_{n=0}^{\infty}a_{n+4} x^{n+4} - 4 x^{-3} \sum_{n=0}^{\infty}a_{n+3} x^{n+3} + 6 x^{-2}\sum_{n=0}^{\infty}a_{n+2} x^{n+2}  \\
 - 4 x^{-1} \sum_{n=0}^{\infty}a_{n+1} x^{n+1} + \sum_{n=0}^{\infty}a_{n} x^{n} = \sum_{n=0}^{\infty} x^{n}
\end{array}
\end{equation} 
\\
The generating function is defined as 
\\
\begin{equation}\label{rec8}
f(x)=\sum_{n=0}^{\infty}a_{n} x^{n}
\end{equation}
\\
Before trying to put equation \ref{rec7} in terms of \ref{rec8} it should be noticed that the righthand side of equation \ref{rec6} is the generating function of the geometrical series whose corresponding succesion is $1,1,1,1,\ldots$. If we multiply each side of equation \ref{rec7} by $x^{4}$ we obtain
\\
\begin{equation}\label{rec9}
\begin{array}{l}
 \sum_{n=0}^{\infty}a_{n+4} x^{n+4} - 4 x \sum_{n=0}^{\infty}a_{n+3} x^{n+3} + 6 x^{2}\sum_{n=0}^{\infty}a_{n+2} x^{n+2}  \\
 - 4 x^{3} \sum_{n=0}^{\infty}a_{n+1} x^{n+1} + x^{4} \sum_{n=0}^{\infty}a_{n} x^{n} = x^{4} \sum_{n=0}^{\infty} x^{n}
\end{array}
\end{equation} 
\\
It can be noticed that by this operation the righthand side of \ref{rec9} correspond to a right shift of the corresponding succession what gives as result the cancelation of the four first places of this succession and the result is the succession $0,0,0,0,1,1,1,1,\ldots$. Rewriting \ref{rec9} in terms of the generating function \ref{rec8} we get
\\
\begin{equation}\label{rec10}
\begin{array}{l}
(f(x)-a_{0}-a_{1}x-a_{2}x^{2}-a_{3}x^{3})-4x(f(x)-a_{0}-a_{1}x-a_{2}x^{2}) \\
+6x^{2}(f(x)-a_{0}-a_{1}x)-4x^{3}(f(x)-a_{0}) + x^{4}f(x)=\frac{x^{4}}{(1-x)} 
\end{array}
\end{equation}
\\
Replacing the $a_{0}=1,a_{1}=2$ and $a_{2}=4,a_{3}=8$ \ref{rec10} we get
\\
\begin{equation}\label{rec11}
\begin{array}{l}
(f(x)-1-2x-4x^{2}-8x^{3})-4x(f(x)-1-2x-4x^{2}) \\
+6x^{2}(f(x)-1-2x)-4x^{3}(f(x)-1) + x^{4}f(x)=\frac{x^{4}}{(1-x)} 
\end{array}
\end{equation}
\\
By algebraic simplification and factorization of \ref{rec11} we obtain
\\
\begin{equation}\label{rec12}
f(x)(1-x)^{4}+(-1+2x-2x^{2})=\frac{x^{4}}{(1-x)}
\end{equation} 
\\
From \ref{rec12} we can obtain $f(x)$ 
\\
\begin{equation}\label{rec13}
f(x)=\frac{x^{4}}{(1-x)^{5}}+\frac{1}{(1-x)^{4}}-\frac{2x}{(1-x)^{4}}+\frac{2x^{2}}{(1-x)^{4}}
\end{equation}
\\ 
By partial fraction decomposition \cite{Grimaldi1} \cite{Sedgewick1} we have that
\\
\begin{equation}\label{partial1}
\frac{x^{4}}{(1-x)^{5}}=\frac{1}{(1-x)}-\frac{4}{(1-x)^{2}}+\frac{6}{(1-x)^{3}}-\frac{4}{(1-x)^{4}}+\frac{1}{(1-x)^{5}}
\end{equation}
\\
\begin{equation}\label{partial2}
\frac{-2x}{(1-x)^{4}}=\frac{2}{(1-x)^{3}}-\frac{2}{(1-x)^{4}}
\end{equation}
\\
\begin{equation}\label{partial3}
\frac{-2x^{2}}{(1-x)^{4}}=\frac{2}{(1-x)^{2}}-\frac{4}{(1-x)^{3}}+\frac{2}{(1-x)^{4}}
\end{equation}
\newpage
Substituing \ref{partial1},\ref{partial2} and \ref{partial3} in \ref{rec13} we get

\begin{equation}\label{rec14}
f(x)=\frac{1}{(1-x)^{5}}-\frac{3}{(1-x)^{4}}+\frac{4}{(1-x)^{3}}-\frac{2}{(1-x)^{2}}+\frac{1}{(1-x)}
\end{equation}
\\
From the generating functions theory it is known that \cite{Grimaldi1} \cite{Sedgewick1}
\\
\begin{equation}\label{rec15}
\begin{array}{l}
\frac{1}{(1-x)^{r}}=\binom{-r}{0}+\binom{-r}{1}(-x)+\binom{-r}{2}(-x)^{2}+\ldots \\
=\sum_{i=0}^{\infty}\binom{-r}{i}x^{i}=1+(-1)\binom{r+1-1}{1}+(-1)^{2}\binom{r+2-1}{2}+\ldots \\
=\sum_{i=0}^{\infty}\binom{r+i-1}{i}x^{i}
\end{array}
\end{equation}
\\
Taking into account \ref{rec15} and using it in the equation \ref{rec14} we can calculate the n-th coeficient of each term and obtain
\\
\begin{equation}\label{rec16}
f(n)=\binom{n+4}{n}- 3 \binom{n+3}{n} + 4 \binom{n+2}{n} - 2 \binom{n+1}{n} + \binom{n}{n} 
\end{equation}
\\
We calculate the polynomials in $n$ from the terms of equation \ref{rec16}
\\
\begin{equation}\label{rec17}
\binom{n}{n}=1
\end{equation}
\\
\begin{equation}\label{rec18}
\binom{n+1}{n}=\frac{(n+1)n!}{n!(n+1-n)!}=n+1
\end{equation}
\\
\begin{equation}\label{rec19}
\binom{n+2}{n}=\frac{(n+2)(n+1)n!}{n!(n+2-n)!}=\frac{n^{2}+3n+1}{2}
\end{equation}
\\
\begin{equation}\label{rec20}
\binom{n+3}{n}=\frac{(n+3)(n+2)(n+1)n!}{n!(n+3-n)!}=\frac{n^{3}+6n^{2}+11n+6}{6}
\end{equation}
\\
\begin{equation}\label{rec21}
\binom{n+4}{n}=\frac{(n+4)(n+3)(n+2)(n+1)n!}{n!(n+4-n)!}=\frac{n^{4}+10n^{3}+35n^{2}+50n+24}{24}
\end{equation}
\\
Replacing \ref{rec17},\ref{rec18},\ref{rec19},\ref{rec20} and \ref{rec21} in equation \ref{rec16} and simplifying we obtain
\\
\begin{equation}\label{rec22}
f(n)=\frac{n^{4}-2n^{3}+11n^{2}+14n-48}{24}
\end{equation}
\\
Recalling that the number of regions $m$ formed in the circle is related with $n$ by $n=m-1$ as well as the fact that for setting the recurrence we have introduced a right shift of four places in \ref{rec9} we replace $n$ by $m-1$ in \ref{rec22} and add $3$ to this new expression for the recuperation of the terms that we have missed because of this right shift in the succession we get
\begin{equation}\label{rec23}
f(m)=\frac{m^{4}-6m^{3}+23m^{2}-18m+24}{24}
\end{equation}

That is the desired result.

\end{proof}
\end{theorem}

\begin{theorem}\label{theoremnonhomogeneous}
Let $n$ the subindex of the recurrence \ref{rec3}, $m$ the number of points on the circle and $f(m)$ the number of regions formed by the division of the circle by chords. Let $n=m-1$. If we solve \ref{rec3} by the solution of linear non-homogeneous with constant coeficients method the number of regions formed by the division of the circle by chords
is
\begin{equation}
f(m) = \frac{m^4-6m^3+23m^2-18m+24}{24} 
\end{equation}

\begin{proof}
From \ref{lemmarecurrence} we have \ref{rec3} .
From the discrete mathematics methods, it is known that the general solution of equation like \ref{rec3} consist of two parts, the solution of the associated  homogeneous equation and a particular solution of \ref{rec3} that is \cite{Grimaldi1} \cite{Sedgewick1}
 
\begin{equation} \label{rec24}
a^{g}_{n}=a^{p}_{n} + a^{h}_{n}
\end{equation}

where the superscript $g$ means general, $p$ means particular and $h$ means homogeneous solutions respectively.
The form of the linear non-homogeneus is

\begin{equation} \label{nh1}
a_{n+4}-4 a_{n+3} + 6 a_{n+2} - 4 a_{n+1} + a_{n} = f(n)
\end{equation}

where $f(n)=1$. Then our non homogeneus equation with the corresponding four initial conditions will be

\begin{equation} \label{nh2}
a_{n+4}-4 a_{n+3} + 6 a_{n+2} - 4 a_{n+1} + a_{n} = 1, ~~a_{0}=1,a_{1}=2,a_{3}=4,a_{4}=8
\end{equation}

Let me start with the particular solution.
As will be shown later, when the homogeneous solution is obtained, the characteristic polynomial have as root $r=1$ with multiplicity four. It means that, in order to have the corresponding 
four linearly independent solutions for the homogeneous associateed equations we will have a polynomial function of third degree as solution \cite{Grimaldi1}. In order to have a particular solution beeing linearly independent from the associated homogeneus solution, the form of the particular solution will be $a^{p}_{n}= An^{4}$ for $A$ constant. Replacing such solution in \ref{nh2} we get

\begin{equation} \label{nh3}
An^{4}-4A(n-1)^4+6A(n-2)^{4}-4A(n-3)^{4}+A(n-4)^4=1
\end{equation}

algebraically developing the fourth degree binomials in \ref{nh3} we obtain
\small{
\begin{equation}\label{nh4}
\begin{array}{l}
A(n^4-4(n^4-4n^3+6n^2-4n+1)+6(n^4-8n^3+24n^2-32n+16) \\
-4(n^4-12n^3+54n^2-108n+1)+(n^4-16n^3+96n^2-256n+256))=1
\end{array}
\end{equation}   
 }

Simplifying \ref{nh4} we get the value
\begin{equation}\label{nh5}
A=\frac{1}{24}
\end{equation}

and from \ref{nh5} we subtitute this value in the particular solution and obtain 

\begin{equation}\label{nh6}
a^{p}_{n}= \frac{n^{4}}{24}
\end{equation}

The associated homogenous equation is

\begin{equation} \label{nh7}
a_{n+4}-4 a_{n+3} + 6 a_{n+2} - 4 a_{n+1} + a_{n} = 0, ~~a_{0}=1,a_{1}=2,a_{3}=4,a_{4}=8
\end{equation}

It is known from the discrete mathematics metodology \cite{Grimaldi1} \cite{Sedgewick1} for solving homogeneous recurrence as \ref{nh7} that their solutions have the form

\begin{equation} \label{nh8}
a^{h}_{n}=Cr^n
\end{equation}

Applying \ref{nh8} to \ref{nh7} it is obtained

\begin{equation} \label{nh9}
Cr^{n+4}-4Cr^{n+3}+6Cr^{n+2}-4Cr^{n+1}+4Cr^{n}=0
\end{equation}

Dividing each member of \ref{nh9} by $Cr^{n}$ we obtain the following characteristic polynomial

\begin{equation} \label{nh10}
r^{4}-4r^{3}+6r^{2}-4r+4=0
\end{equation}

By factorization of \ref{nh10} we get

\begin{equation} \label{nh11}
(r-1)^4=0
\end{equation}

Then the roots of \ref{nh11} are $r=1$ with multiplicity of four. In order to have four linearly independent homogeneous solutions the homogeneous solution will have the next form

\begin{equation} \label{nh12}
a^{h}_{n}=C_{1}1^{n}+C_{2}n1^{n}+C_{3}n^{2}1^{n}+C_{4}n^{3}1^{n}=C_{1}+C_{2}n+C_{3}n^{2}+C_{4}n^{3}
\end{equation} 

Where $C_{1},C_{2},C_{3}$ and $C_{4}$ are constants that can be determined by the application of the
initial conditions $a_{0}=1,a_{1}=2,a_{3}=4$ and $a_{4}=8$. From \ref{rec24} we know that $a^{g}_{n}=a^{p}_{n}+a^{h}_{n}$ and the we can establis the following relation

\begin{equation} \label{nh13}
a^{g}_{n}=a^{h}_{n}+a^{p}_{n}=C_{1}+C_{2}n+C_{3}n^{2}+C_{4}n^{3}+\frac{n^4}{24}
\end{equation} 

By application of the initial conditions to \ref{nh13} we establish the following relations
\begin{equation} \label{nh14}
\begin{array}{l}
a_{0}=C_{1}+C_{2}(0)+C_{3}(0)^2+C_{4}(0)^3+\frac{(0)^4}{24}=1\\
a_{1}=C_{1}+C_{2}(1)+C_{3}(1)^2+C_{4}(1)^3+\frac{(1)^4}{24}=2\\
a_{2}=C_{1}+C_{2}(2)+C_{3}(2)^2+C_{4}(2)^3+\frac{(2)^4}{24}=4\\
a_{3}=C_{1}+C_{2}(3)+C_{3}(3)^2+C_{4}(3)^3+\frac{(3)^4}{24}=8
\end{array}
\end{equation}

From the relations \ref{nh14} we establish the following linear relations

\begin{equation} \label{nh15}
\begin{array}{l}
C_{1}=1\\
C_{2}+C_{3}+C_{4}=\frac{23}{24}\\
2C_{2}+4C_{3}+8C_{4}=\frac{56}{24}\\
3C_{2}+9C_{3}+27C_{4}=\frac{87}{24}
\end{array}
\end{equation}

For obtaining the values of $C_{2},C_{3}$ and $C_{4}$ we solve the following linear system

\begin{equation}\label{nh16}
\left( \begin{array}{ccc}
1 & 1 & 1 \\
2 & 4 & 8 \\
3 & 9 & 27 \\
\end{array}
\right)
\left(
\begin{array}{c}
C_{2}\\
C_{3}\\
C_{4} 
\end{array}
\right)
=
\left(
\begin{array}{c}
\frac{23}{24}\\
\frac{56}{24}\\
\frac{87}{24} 
\end{array}
\right)
\end{equation}

Solving the linear system \ref{nh16} we obtain

\begin{equation}\label{nh17}
\begin{array}{c}
C_{2}=\frac{14}{24}\\
C_{3}=\frac{11}{24}\\
C_{4}=\frac{-2}{24} 
\end{array}
\end{equation}

Replacing the values obtained in \ref{nh13} we get the following expression

\begin{equation} \label{nh18}
a^{g}_{n}=1+\frac{14}{24}n+\frac{11}{24}n^{2}-\frac{2}{24}n^{3}+\frac{n^4}{24}
\end{equation}

We know that the relation between the number $m$ of points in the circle from where the chords are traced, and
$n$ is $n=m-1$ so we can express \ref{nh18} in terms of $m$ as follows

\begin{equation}\label{nh19}
a^{g}_{n}=1+\frac{14}{24}(m-1)+\frac{11}{24}(m-1)^{2}-\frac{2}{24}(m-1)^{3}+\frac{(m-1)^4}{24}
\end{equation}

Algebraically developing each term of \ref{nh19} we get

\begin{equation}\label{nh20}
\begin{array}{l}
a^{g}_{n}=1+\frac{14}{24}m-\frac{14}{24}+\frac{11}{24}(m^{2}-2m+1)\\
-\frac{2}{24}(m^{3}-3m^{2}+3m-1)+\frac{(m^4-4m^3+6m^2-4m+1)}{24}
\end{array}
\end{equation}

Simplifying \ref{nh20} we finally obtain

\begin{equation} \label{nh21}
a^{g}_{n}=\frac{24-18m+23m^2-6m^3+m^4}{24}
\end{equation}

Reversing the order of the numerator terms we finally obtain the desired result

\begin{equation} \label{nh22}
f(m)=\frac{m^4-6m^3+23m^2-18m+24}{24}
\end{equation}

\end{proof}
\end{theorem}

\section{Conclusions}

In this article, I proposed a deduction of the Leo Moser's known formula for counting the number of regions that are formed by dividing the circle by chords, by solving a fourth order difference equation obtained by the successive differences method. I solved this recurrence equation by two different methods. This article illustrate how a classical problem can lead to different and creative developments. That is what makes mathematics so exciting.  



\begin{thebibliography}{1}






\bibitem{Grimaldi1}
{\small {\sc Ralph P. Grimaldi:} {\it Discrete and Combinatorial Mathematics: An applied introduction}
Addison-Wesley, {\bf 3th Ed.} (1994), .}

\bibitem{Jobbings1}
{\small {\sc Andrew Jobbings:} {\it A selection of mathematical articles and notes by Andrew Jobbings.}
,http://www.arbelos.co.uk/Papers/Chords-regions.pdf,{\bf 27 December} (2008)}

\bibitem{Maier1}
{\small {\sc Eugene Maier:} {\it Counting Pizza Pieces and Other Combinatorial Problems.}
Mathematics Teacher, {\bf 81} (1988), 22--26.}

\bibitem{Moser1}
{\small {\sc Leo Moser and W. Bruce Ross:} {\it Mathematical Miscellany.}
Mathematics Magazine, {\bf 23} (1949), 109--114.}

\bibitem{Sedgewick1}
{\small {\sc Robert Segdewick and Philippe Flajolet:} {\it Introduction to the Analysis of Algorithms.}
Addison-Wesley, {\bf 2nd Printing} (2001), .}

\end{thebibliography}

\end{document}